
%
\input amstex
\documentstyle{amsppt}
\NoBlackBoxes
\magnification = \magstep1

\nologo

\NoRunningHeads

\TagsOnRight

\pageheight{7.4in} \pagewidth{6.5truein}

\topmatter
\title
Variation of argument and Bernstein index for holomorphic
functions on Riemann surfaces
\endtitle

\rightheadtext{}

\author      Yu. Ilyashenko
\endauthor

\address
Cornell University, US; Moscow State and Independent Universities,
Steklov Math. Institute, Moscow
\endaddress

\thanks
The author was supported by part by the grants NSF 0400495, RFBR
02-02-00482.
\endthanks

\abstract An upper bound of the variation of argument of a
holomorphic function along a curve on a Riemann surface is given.
This bound is expressed through the Bernstein index of the
function multiplied by a geometric constant. The Bernstein index
characterizes growth of the function from a smaller domain to a
larger one. The geometric constant in the estimate is explicitly
given. This result is applied in \cite {GI} to the solution of the
restricted version of the infinitesimal Hilbert 16th problem,
namely, to upper estimates of the number of zeros of abelian
integrals in complex domains.
\endabstract

\endtopmatter

\define\e {\varepsilon }
\define\nbd {neighborhood }
\define\ph {\varphi }

\document

\head Introduction \endhead

Consider a holomorphic function $f$ in a topological disk $U$
imbedded in a Riemann surface, and a curve $\Gamma \subset U.$ We
are interested in an upper bound of the variation of argument of
the function $f$ along the curve $\Gamma $. This bound may be
expressed through the growth rate of $f$ from some intermediate
domain $U''$ to $U,$ and through a geometric factor depending on
the relative position of all the three sets:  $\Gamma \subset U''
\subset U.$ Our result improves the theorem from \cite {KYa} that
provides a similar estimate for a domain $U \subset \Bbb C$, and
without explicit formula for the geometric factor mentioned above.
We provide a formula for this factor, and consider functions on a
Riemann surface. Both improvements are crucial for an application
of our main result to an upper estimate of the number of zeroes of
Abelian integrals in the complex domain. This latter estimate
forms a part of the solution of so called {\it restricted
infinitesimal Hilbert 16th problem,} \cite {GI}.

Our main tool is Growth-and-Zeros theorem stated and improved in
the next section. The main results of the paper are stated in
Section 2.

\head \S 1. Bernstein index and Growth-and-Zeros theorem
\endhead

 Let $W$ be a Riemann surface, $\pi : W \to \Bbb
C$ be a holomorphic function (called projection) with non-zero
derivative. Let $\rho$ be the metric on $W$ lifted from $\Bbb C$
by projection $\pi .$ Let $U \subset W$ be a connected domain, and
$K \subset U $ be a compact set. For any $p \in U$ let $\e (p,
\partial U)$ be the supremum of radii of the disks centered at $p,$
located in $U$ and such that $\pi $ is bijective on these disks.
{\it The $\pi $-gap between $K$ and $\partial U,$} is defined as
$$
\pi \text {-gap } (K, \partial U) = \min_{p\in K} \e (p,\partial
U).
$$
\

\proclaim {Growth-and-zeros theorem } Let $W, \pi , \rho$ be the
same as before. Let $U \subset W$ be a domain conformally
equivalent to a disk. Let $K \subset U$ be a path connected
compact subset of $U$ (different from a single point). Suppose
that the following two assumptions hold:

Diameter condition:
$$
{\text { diam }}_{int} K \le D;
$$

Gap condition:
$$
\pi \text {-gap} (K, \partial U) \le \e .
$$

Let $I$ be a bounded holomorphic function on $\bar U.$ Then
$$
\# \{ z \in K|I(z) = 0\} \le e^{\frac {2D}{\e }}\log\frac
{\max_U|I|}{\max_K|I|}   \tag 1
$$
\endproclaim
The definition  of the intrinsic diameter is well known; yet we
recall it for the sake of completeness.
\definition {Definition 1} The {\it intrinsic distance} between two
points of a path connected set in a metric space is the infinum of
the length of paths in $K$ that connect these points (if exists).
The {\it intrinsic diameter} of $K$ is a supremum of intrinsic
distances between two points taken over all the pairs of points in
$K.$
\enddefinition
\definition {Definition 2} The second factor in the right hand side
of (1) is called  {\it the Bernstein index } of $I$ with respect
to $U$ and $K$ and denoted by $B_{K,U}(I):$
$$
B_{K,U}(I) = \log \frac {M}{m}, \ M = \sup_U|I|, \ m = \max_K|I|.
$$
\enddefinition

\demo {Proof of the Growth-and-Zeros theorem} The above theorem is
proved in \cite {IYa} for the case when $W = \Bbb C, \pi = id.$ In
fact, in \cite {I99} another version of (1) is proved with (1)
replaced by
$$
\# \{ z \in K|I(z) = 0\} \le B_{K,U}(I)e^\rho ,  \tag 2
$$
where $\rho $ is the diameter of $K$ in the Poincar\' e metric of
$U.$ In this case it does not matter whether $U$ belongs to $\Bbb
C$ or to the Riemann surface.

\proclaim {Proposition 1} Let $K, U$ be two sets in the Riemann
surface $W$ from Definition 2, and let the Diameter and Gap
conditions from the Growth-and-Zeros theorem hold. Then the
diameter of $K$ in the Poincar\' e metric of $U$ admits the
following upper estimate:

$$
\rho \le \frac {2D}{\e} .   \tag 3
$$
\endproclaim

Growth-and-Zeros theorem now follows from (2) and Proposition 1.
\enddemo

\demo {Proof of Proposition 1} By the gap condition, for any point
$p \in K,$ the domain $U$ contains a disk $D$ of radius $\e$
centered at $p.$ In more detail, $D \subset U $ is a set mapped by
$\pi$ bijectively onto a disk of radius $\e$ in $\Bbb C$ centered
at $\pi (p)$ and such that $p \in D.$ {\it In what follows, for
any topological disk $V$ the Poincar\'e matric of $V$ is denoted
by $PV.$} The  length of a vector $v$ attached at $p$ in the
metric $PD$ equals to the Eucleadian length $|v|$ of $v$ divided
by $\e$ and multiplied by 2. By the monotonicity property of the
Poincar\'e metric, the length of $v$ in sense of $PU$ is no
greater then the previous one. This implies (3).
\enddemo

\head \S 2. Upper estimate of the variation of argument
\endhead

A definition of a {\it variation of argument of a complex valued
function along a curve } is contained in this, rather lengthy,
name; this variation is denoted $V_{\Gamma}(f)$ for a function $f$
and a curve $\Gamma.$ In more detail, let $U$ be a domain on a
Riemann surface,  $\Gamma : [0,1] \to U $ be a curve, $f: U \to
\Bbb C ,$ $f|_{\Gamma} \not = 0$  a function. Fix an arbitrary
branch of $\arg f$ on $\Gamma $ and let $ \phi = \arg f \circ
\Gamma.$ By definition, $V_{\Gamma}(f)$ is the total variation of
$\phi $ on [0,1].

A first step in  establishing a relation between variation of
argument and the Benstein index is done by the following KYa
(Khovansky-Yakovenko) theorem. Let $U, V $ and $f$ be the same as
in the previous paragraph.

\proclaim {KYa theorem, \cite {KYa}} For any tuple $U,\Gamma
\subset U $ as above and a  compact set $K\subset U$ there exists
a geometric constant $\alpha=\alpha (U,K,\Gamma),$ such that
$$ V_{\Gamma}(f)\le \alpha B_{K,U}(f).$$
\endproclaim

In  \cite {KYa} an upper estimate of the    Bernstein index
through the variation of the argument along $ \Gamma = \partial U$
is given; we do not use this estimate. On the contrary, we need an
improved version of the previous theorem with $\alpha$ explicitly
written and $ U$ being a domain on  a Riemann surface. These two
improvements  are achieved in the following two theorems.

Let $|\Gamma|$ be the length, and $ \kappa (\Gamma ) $ be the
total curvature of a curve on a surface endowed with a Riemann
metric.

\proclaim {Theorem 1} Let $\Gamma \subset U'' \subset U' \subset U
\subset \Bbb C$ be respectively a piecewise smooth curve, and
three path connected  open sets in $\Bbb C, \  \bar U $ is compact
and simply connected. Let $f: \bar U \to \Bbb C$ be a holomorphic
function, $f|_{\Gamma} \not = 0.$ Let $ \frac {D}{\e}
> 3,$ and the following two conditions hold:

Gap condition:
$$
\rho (\Gamma , \partial U'') \ge \e , \ \rho (U'', \partial U')
\ge \e , \ \rho (U',\partial U) \ge \e; \tag 4
$$

Diameter condition:

$$
\text { diam }_{\text {int}}U'' \le D, \ \text { diam }_{\text
{int}}U' \le D.  \tag 5
$$
Then
$$
V_\Gamma (f) \le B_{U'',U}(f)(\frac {\mid \Gamma \mid}{\e} +
\kappa (\Gamma ) + 1)e^{\frac {5D}{\e }}.   \tag 6
$$
\endproclaim
\proclaim {Theorem 2} Let $\Gamma \subset U'' \subset U' \subset U
\subset W$ be respectively a piecewise smooth curve, and three
path connected open sets in a Riemann surface $W, \ \bar U $ is
compact  and  simply connected. Let $f: \bar U \to \Bbb C$ be a
holomorphic function, $f|_{\Gamma} \not = 0.$ Let $\pi: W \to \Bbb
C$ be a projection which is locally biholomorphic, and the metric
on $W$ be a pullback of that on $\Bbb C.$ Let $ \frac {D}{\e}
> 3,$ and the following two conditions hold:

Gap condition:
$$
\pi \text {-gap }(\Gamma , U'') \ge \e , \ \pi \text {-gap
}(U'',U') \ge \e , \ \pi \text {-gap }(U',U) \ge \e ;   \tag 7
$$

Diameter condition:

$$
{\text {diam }}_{int}U''\le D, \ {\text {diam }}_{int}U'\le D.
\tag 8
$$
Then inequality (6) holds.
\endproclaim

The proofs of these theorems is partly based on the methods of
\cite {KYa}. Yet our presentation is self contained. Theorem 1
improves the parallel result of \cite {KYa} by giving an explicit
expression for the geometric factor in the KYa Theorem. Theorem 2
has no analoge in \cite {KYa}. Note that in \cite {KYa} the
variation of argument (modulo a factor $( 2 \pi)^{-1}$) is called
a {\it Voorhoeve index.}

\head \S 3. Variation of argument of nowhere zero functions on
Riemann surfaces
\endhead

\proclaim {Lemma 1} Let $\Gamma ,U'',U' \subset W$ be the same as
in Theorem 2, and $F: U' \to \Bbb C$ be a nowhere zero holomorphic
function. Then
$$
V_\Gamma (F) \le B_{U'',U'}(F) \frac {|\Gamma|}{\e} e^{\frac
{2D}{\e }}. \tag 9
$$
\endproclaim
\demo {Proof} The classical definition of the variation of a
function of one variable implies:
$$
V_{\Gamma}(f)=\int_{\Gamma}|Im{{(\log f)}^{\prime}}|ds,
$$
$s$ is the arc length parameter on $\Gamma$. Hence,
$$
V_\Gamma (F) \le \int_\Gamma \mid (\log F)'\mid ds.    \tag 10
$$
Without loss of generality, we may assume that
$$
\max_{\bar U''} \mid F \mid = 1,
$$
and $F(a) = 1$ for some $a \in U''.$ Then, for $B' =
B_{U'',U'}(F)$ we have: $(\log F)(U') \subset \Bbb C_{B'}^-:= \{ w
\mid \text { Re }w \le B'\} .$ Let $g = \frac {1}{B'}\log F.$
Then, by (10),
$$
V_\Gamma (F) \le \mid \Gamma \mid \cdot B' \cdot \max_\Gamma \mid
g'\mid .  \tag 11
$$
Note that
$$
g(U') \subset \Bbb C_1^-.  \tag 12
$$

We want  to get an upper bound for $\mid g' \mid $ on $\Gamma$
making use of the Cauchy estimates. For this we need to estimate
from above $\max_{z\in U''}\mid g(z)\mid .$ To do that we will use
an upper  estimate  of the diameter of $U''$ in the Poincar\' e
metric $PU'$ of $U',$ and the fact that $g$ does not increase the
Poincar\' e metric. By Proposition 1, we have:
$$
{\text { diam }}_{PU'}U'' \le \frac {2D}{\e }.
$$
Note that for the above chosen $a \in U'', \ g(a) = 0.$ Then, by
(12) and monotonicity property of the Poincar\' e metric, for any
$z \in U'',$
$$
\rho_{P\Bbb C_1^-}(0,g(z)) \le \rho_{PU'}(a,z) \le \frac {2D}{\e
}.
$$
This estimate for the Poincar\' e distance between $0$ and $g(z)$
in sense of $P\Bbb C_1^-$ implies an estimate on the Euclidean
distance between $0$ and $g(z):$
$$
\mid g(z)\mid \le e^{\frac {2D}{\e }}.
$$
By the Cauchy estimate and gap condition (4) for $\Gamma $ and
$U'',$
$$
\max_\Gamma \mid g'\mid \le \e ^{-1}e^{\frac {2D}{\e }}.
$$
Together with (11), this implies the lemma.
\enddemo

\head  \S 4.  Variation of argument of holomorphic functions on
$\Bbb C $ \endhead

\demo {Proof of the Theorem 1} Lemma 1 implies the theorem in the
case when $f \ne 0$ in $U'.$ The lemma should be applied to $F =
f\mid _{U'}$ keeping in mind that $ B_{U'', U'} (F) \le B_{U'', U}
(F).$ In  general, $f$ may have zeros in $U'.$ Let $d = \# \{ f(z)
= 0\mid z \in U'\},$ zeros of $f$ are counted with multiplicities.
 By the Growth-and-Zeros theorem,
$$
d \le B_{U',U}(f)e^{\frac {2D}{\e }}.
$$
 Let $ B =
B_{U'',U}(f).$ Note that $B_{U',U}(f) \le B.$ Hence,
$$
d \le B e^{\frac {2D}{\e }}. \tag 13
$$
Let $p$ be a monic polynomial of degree $d$ that has the same
zeros as $f$ in $U',$ with the same multiplicities. Then the
function
$$
F = \frac {f}{p}
$$
is holomorphic and nowhere zero in $U'.$ Equality $f = Fp$ implies
$$
V_\Gamma (f) \le V_\Gamma (F) + V_\Gamma (p).    \tag 14
$$
The proof is based on the triangle inequality applied to moduli of
derivatives  of $\arg F$ and $\arg p$ along $\Gamma.$ In more
detail, (14) is proved in \cite {KYa}.

\proclaim {Lemma 2} Let $U''\subset U' \subset U \subset \Bbb C$
be the same as in Theorem 1. Let $F$ be a holomorphic function in
$\bar U, \ p$ be a monic polynomial of degree $d$ with zeros
located in $U'.$ Let $f = Fp, \ B = B_{U'',U}(f),$ as before. Then
$$
B_{U'',U}(F) \le B + d\log \frac {D}{\e }.   \tag 15
$$
\endproclaim

This lemma is very close to a more complicated statement proved in
\cite {KYa}. The proof of the lemma itself follows; it is
straightforward. Before the proof, note that the choice of the
domains in Lemma 2 is important. On one hand,  the same Bernstein
index $B$ is used in (13) and (15). On the other hand, the right
hand side of (15) estimates from above the Bernstein index $B' $
used in Lemma 1.

\demo {Proof of Lemma 2} Without loss of generality we may assume
that
$$
\max_{\bar U''} |f| = 1.
$$
Then
$$
\log \max_{\bar U} |f| = B.
$$
Hence,
$$
\log \max_{\partial U}|F| \le B - \log \min_{\partial U} |p|.
$$
On the other hand, $p$ is a product of $d$ binomials of the form
$w - w_j$ with $w_j \in U'.$ Then, by the gap condition (4) for
$U'$ and $U,$
$$
 \log \min_{\partial U} |p| \ge - d\log \e .
$$
Hence,
$$
\log \max_{\bar U} |F| \le B - d \log \e
$$
by the maximum modulus principle.

On $U'', \ |p| \le D^d.$ Then $\log \max_{U''}|F| \ge - d\log D.$
Hence,
$$
B_{U'',U}(F) \le B + d(\log D - \log \e ).
$$
This proves the lemma.
\enddemo

Now let us complete the proof of Theorem 1.  To do that, let us
estimate the terms in the right hand side of (14). First, by \cite
{KYa},
$$
V_\Gamma (p) \le (\kappa (\Gamma ) + 2\pi )d.
$$
The proof of this inequality is based on the decomposition of a
polynomial into a product of linear factors. For $q$ linear, an
elementary estimate $V_\Gamma (q) \le \kappa (\Gamma ) + 2 \pi$
holds. By (13),
$$
V_\Gamma (p) \le (\kappa (\Gamma ) + 2\pi )B e^{\frac {2D}{\e }}
\tag 16
$$

Second,  $ V_\Gamma (F)$ is estimated in Lemma 1, see (9). But
$$ B_{U'', U'} (F) \le B_{U'', U} (F).
$$
The  Bernstein index in the right hand side is already estimated
from above in (15). Namely, let $\alpha = \log\frac {D}{\e }.$
Then, by Lemma 2,
$$
B_{U'', U} (F) \le B + \alpha d.
$$
From this, by (9) and (13) we get:
$$
V_\Gamma (F) \le B(1 + \alpha e^{\frac {2D}{\e }})\frac
{|\Gamma|}{\e}e^{\frac {2D}{\e }}.
$$

By assumption of Theorem 1, $ \frac {D}{\e}
> 3. $  Elementary estimates yield:
$$
1 + \alpha e^{\frac {2D}{\e }} \le  e^{\frac {3D}{\e }}.
$$
Hence,
$$
V_\Gamma (F) \le B \frac {\mid \Gamma \mid}{\e} e^{\frac {5D}{\e
}}. \tag 17
$$
By (14), (16) and (17) we get now
$$
V_\Gamma (f) \le B( \frac {\mid \Gamma \mid}{\e} + \kappa (\Gamma
) + 2\pi )e^{\frac {5D}{\e }}.
$$
This proves Theorem 1.
\enddemo

\head \S 5. Variation of argument of holomorphic functions on
Riemann surfaces \endhead

\demo {Proof of Theorem 2} Theorem 2 is proved in this and the
next sections.

Let $\ph : U \to D_1$ be a conformal mapping of $U$ onto the unit
disk $D_1 $ that takes zero value at some point $b \in U''.$ Let
$$
g = f\circ \ph^{-1}: D_1 \to \Bbb C.
$$
The number of zeros, the Bernstein index and the variation of the
argument of a holomorphic function are invariant under
biholomorphic maps of the domain of a function. In particular, the
function $g$ has the same number $d$ of zeros, counted with
multiplicities, in $\ph (U'),$ as $f$ has in $U'.$ Note that
estimate (13) for $d$ holds by Growth-and-Zeros theorem because
the theorem  is stated for Riemann surfaces.

Let us take a polynomial $p$ that has the same zeros with the same
multiplicities as $g$ in $\ph (U').$  Then $\deg p = d, $ and (13)
holds.  Take $P = p \circ \ph : U \to \Bbb C $ and $F = \frac
{f}{P}.$ Then $F$ is nowhere zero in $U'.$ By the same argument as
above,
$$
V_\Gamma (f) \le V_\Gamma (F) + V_\Gamma (P).   \tag 18
$$

In this section the first term in the right hand side of (18) is
estimated from above. The next section deals with the second term.

The variation  $ V_\Gamma (F) $ is estimated from above by Lemmas
1 and 2, with domains $U'', U', U $ replaced by $ \ph (U''),  \ph
(U'),  \ph (U) = D_1.$ The main part of the proof deals with the
diameter and gap conditions for this new triple of domains.

 By an obvious inequality for Bernstein indexes, and by the invariance of the
Bernstein index mentioned above,
$$
B_{U'',U'}(F)  \le B_{U'',U}(F)= B_{\ph (U''),D_1}\left( \frac
{g}{p}\right).
$$
Once more, by the invariance of the Bernstein index,
$$
B_{\ph (U''), D_1} = B.
$$
Hence, by Lemma 2,
$$
B_{U'',U'}(F)  \le B + d\alpha'
$$
where $\alpha' = \log \frac {D'}{\e'}, \ D' = {{\text
{diam}}_{\text {int}}\ph (U')}, \ \e' = {\text {gap}(\ph
(U'),D_1)}.$ The factor $\alpha' $ is estimated from above in
Proposition 2 below. Hence, by (13),
$$
B_{U'',U'}(F) \le B + \alpha' Be^{\frac {2D}{\e }}.
$$

By Lemma 1
$$
V_\Gamma (F) \le B( 1 + \alpha'e^{\frac {2D}{\e }} )|\Gamma
|\e^{-1}e^{\frac {2D}{\e }}. \tag 19
$$

\proclaim {Proposition  2} In the assumption above,
$$
D' = {\text {diam }}_{\text {int}}\ph (U') \le \frac {D}{\e }
$$
$$
 \e' = \text {gap }(\ph (U'),D_1) \ge e^{-\frac {2D}{\e }}
$$
\endproclaim
\demo {Proof } By the gap condition (7),
$$
\pi \text {-gap } (U',U) \ge \e .
$$
Hence, we may apply Cauchy inequality to $\ph $ and get: for any
$a \in U',$
$$
\mid \ph' (a)\mid \le \e^{-1}.
$$
Here and below the derivatives are taken with respect to the local
parameter $z$ on $W$ lifted by $\pi $ from $\Bbb C.$ Together with
Diameter condition (8), this proves the first statement of the
proposition.

By Proposition 1, $ \text {diam }_{PU} (U') \le \frac {2D}{\e}.$
By the invariance of the Poincar\' e metric under $\ph ,$ we get
that for any $w \in \ph (U'), \ \rho_{PD_1}(0,w) \le \frac {2D}{\e
}.$ Hence,
$$
1 - |w| \ge e^{-\frac {2D}{\e }}.
$$
This proves the  proposition.
\enddemo
\proclaim {Corollary 1} In (19),
$$
\alpha' \le \frac {3D}{\e }. \tag 20
$$
\endproclaim
Finally, by (19) and (20),
$$
V_\Gamma (F) \le  B(1 + \frac {3D}{\e} e^{\frac {2D}{\e }})|\Gamma
|\e^{-1}e^{\frac {2D}{\e }}. \tag 21
$$
\enddemo

\head \S 6. Total curvature and conformal mappings \endhead

Let us now estimate the second term in (18). By the invariance of
the variation of argument,
$$
V_\Gamma (P) = V_{\ph (\Gamma )}(p) \le (\kappa (\ph (\Gamma )) +
2\pi )d.   \tag 22
$$
Once more, the only thing to do is to estimate from above the
total curvature of the new curve $ \ph (\Gamma ).$

\proclaim {Lemma 3} In the assumption above,
$$
\kappa (\ph (\Gamma )) \le \kappa (\Gamma ) + 2 \e^{-1}\mid \Gamma
\mid .     \tag 23
$$
\endproclaim
\demo {Proof} We will treat the case when $\Gamma $ is not only
piecewise smooth but smooth. The contribution to $\kappa (\Gamma
)$ given by the vertexes of the piecewise smooth curve remains
unchanged under $\ph ,$ because $\ph $ is a conformal mapping. Let
$\Gamma = \{ \gamma (s)\mid s \in [0,|\Gamma |]\} , \ s$ be a
natural parameter on $\Gamma .$ Let dot denote derivative in $s.$
Then $|\dot \gamma (s)| \equiv 1.$ Let $ s_1 $ be a natural
parameter on $\ph (\Gamma ).$ Then
$$
ds_1 = \mid \ph' \circ \gamma (s)\mid ds.
$$
Let $ \ph \circ \gamma (s) = (x(s), y(s)), \ \kappa_{\Gamma } (s)$
be the curvature of the curve $\Gamma $ at a point $s.$ Then
$$
\kappa_{\ph (\Gamma )}(s) = \frac {|\ddot x\dot y - \dot x\ddot y
|}{{(\dot x^2 + \dot y^2)}^{3/2}}.
$$
Hence,
$$
\kappa ( \ph (\Gamma)) = \int_0^{|\ph (\Gamma )|}\kappa_{\ph
(\Gamma )}(s(s_1))ds_1 = \int_0^{|\Gamma |}\kappa_{\ph (\Gamma
)}(s)\mid \ph' \circ \gamma (s)\mid ds .
$$
On the other hand,
$$
 |\ddot x\dot y - \dot x\ddot y| \le |(\ddot x, \ddot y)| |(\dot
x, \dot y)|, \ {\dot x}^2 + {\dot y}^2 = {|(\dot x, \dot y)|}^2,
$$
$$
\mid (\dot x, \dot y)\mid (s) = \mid \ph' \circ \gamma (s)\mid
\cdot |\dot \gamma (s)| = \mid \ph' \circ \gamma (s)\mid .
$$
$$
|(\ddot x, \ddot y)|(s) = \mid \ph'' \circ \gamma \cdot {\dot
\gamma }^2 + \ph' \circ \gamma \cdot \ddot \gamma \mid (s) \le
\mid \ph'' \circ \gamma \mid (s) + \kappa_\Gamma (s)\cdot |\ph'
\circ \gamma |.
$$
Hence,
$$ \split
\kappa (\ph (\Gamma )) = \int_0^{|\Gamma |}\frac {|\ddot x\dot y -
\dot x\ddot y |}{{\dot x}^2 + {\dot y}^2}ds  \le \int_0^{|\Gamma
|}\left( \frac {\mid (\ph'' \circ \gamma )\mid }{\mid (\ph' \circ
\gamma )\mid }(s) + \kappa_\Gamma (s)\right) ds \\ \ \ \ \ \ \ \ \
\ \ \ \ \ \ \ \ \ \  \ \ \ \ \ \ \ \ \ \ \   \le \max_\Gamma
\left| \frac {\ph'' }{\ph' }\right| \cdot |\Gamma | + \kappa
(\Gamma ).
\endsplit \tag 24
$$
For any $a \in \Gamma ,$ the function
$$
\psi (z) = \frac {\ph (\e z + a)}{\e \ph' (a)}
$$
maps the unit disc conformally into $\Bbb C,$ and $\psi' (0) = 1.$
By the Koebe theorem, \cite {C, \S 8 Ch 6}, we have:  $\mid
\psi''(0)\mid \le 2.$ But
$$
\psi'' (0) = \frac {\ph'' (a)\e }{\ph' (a)}.
$$
Hence,
$$
\left| \frac {\ph'' (a)}{\ph' (a)}\right| \le 2\e^{-1}.
$$
This estimate substituted to (24) implies (23). Lemma 3 is proved.
\enddemo

We can now complete the proof of Theorem 2. By (22) and (13),
$$
V_\Gamma (P) \le B(\kappa (\ph (\Gamma )) + 2\pi )e^{\frac {2D}{\e
}}.
$$
By Lemma 3,
$$
V_\Gamma (P) \le B(\kappa (\Gamma ) + 2\frac {|\Gamma |}{\e} +
2\pi ) e^{\frac {2D}{\e }}.
$$
By (18), (21), assumption $ D/\e > 3,$  and the previous
inequality,
$$
V_\Gamma (f) \le Be^{\frac {2D}{\e }}\left[ 2\pi + \kappa (\Gamma
) + |\Gamma | \left( \frac {3}{\e } + \frac {3D}{\e^2}e^{\frac
{2D}{\e }}\right) \right] \le Be^{\frac {5D}{\e }}\left (1 +
\kappa (\Gamma ) + \frac {|\Gamma |}{\e}\right).
$$
This proves Theorem 2.

\enddocument